\documentclass[12pt]{amsart}[1996/10/24]
\relax
\usepackage{amsmath,amscd,stmaryrd}

\advance \topmargin by -\headheight
\advance \topmargin by -\headsep
\evensidemargin \oddsidemargin
\newtheorem{theorem}[equation]{Theorem}
\newtheorem{prop}[equation]{Proposition}
\newtheorem{lemma}[equation]{Lemma}
\newtheorem{corollary}[equation]{Corollary}
\newtheorem{definition}[equation]{Definition}

\numberwithin{equation}{section}
\newcommand{\Vol}{\operatorname{Vol}}
\newcommand{\Rm}{\operatorname{Rm}}

\newcommand{\Ric}{\operatorname{Ric}}
\newcommand{\Diam}{\operatorname{Diam}}
\newcommand{\Hess}{\operatorname{Hess}}

\newcommand{\be}{\begin{equation}}
\newcommand{\ee}{\end{equation}}

\begin{document}
\title[]
{A remark on degenerate singularity in three dimensional Ricci flow}
\author{Yu Ding}
\address{Department of Mathematics and Statistics, 
California State University, Long Beach, CA 90840}

\email{yding@csulb.edu}
\begin{abstract}
We show that a rescale limit at any degenerate singularity of Ricci flow 
in dimension $3$ is a steady gradient soliton. In particular, we give a 
geometric description of type I and type II singularities.             
\end{abstract}
\maketitle
\tableofcontents


Degenerate singularity of the Ricci flow 
\begin{equation}
\frac{d}{dt}g(t)=-2\Ric(g(t)) 
\end{equation}
was introduced in Hamilton's paper 
\cite{HaSing}. In that paper, Hamilton first described the 
{\em nondegenerate neck-pinching}. Roughly speaking, one starts  
the Ricci flow on a dumbbell 
shaped $3$-manifold, with the neck diffeomorphic to $S^2\times [-1, 1]$. 
It is expected that the neck shrinks in the $S^2$ direction,   
where the curvature is very positive, and, at the same time,  
stays relatively stationary in the ${\mathbb R}$ direction, 
where the curvature is slightly negative. 
After some time, the neck pinches off and forms a singularity. 
One step further, Hamilton purposed the notion of {\em degenerate 
neck-pinching }: reduce the left half of the dumbbell  
into a critical size and then start the Ricci flow.  
It is expected that after some time, all of the left half of the dumbbell 
pinches off, 
and forms a singularity like a horn growing out of the 
(remaining) right half of the dumbbell. 
See \cite{HaSing} for further descriptions and 
some very inspiring pictures. 

In this paper, we prove that in dimension $3$, a rescale limit of a  
degenerate singularity of Ricci flow is a steady 
gradient soliton, see Theorem 
\ref{th-DS-main-theorem-steady-soliton-faih-g9uj3-gve89}. 
The precise definition of degenerate singularity, given 
in Definition 
\ref{Define-DS-degenerate-singularity-7834h-8as1}, 
is based on Perelman's notion of {\em canonical neighborhood}. 
Our definition is a geometric one that 
reflects Hamilton's original picture in \cite{HaSing}. On the other hand, 
as we will see later in this paper, this geometric 
definition is equivalent to 
that the singularity being of type II. 

For previous works on neck-pinching, see \cite{AK1}, \cite{AK2}, 
and the book \cite{CLN}. In the book \cite{CK}, 
there is a detailed treatment of nondegenerate neck-pinching 
in chapter 2, and a discussion of degenerate neck pinching 
in page 62-66. 

We start by reviewing some of Perelman's results in \cite{P1}, \cite{P2}; 
for more details, see 
\cite{CaoZhu}, \cite{KL} and \cite{MT}.  
In section 2 we use an estimate on Perelman's $l$ functional to 
rule out noncompact ancient solutions with 
positive curvature that develops a type I singularity. Therefore 
a rescale limit of a degenerate singularity is either 
an eternal solution or an ancient solution that develops a 
type II singularity. In both cases, we need to take a further rescale 
limit in {\em forward time};  
we treat certain issues related to this in Section 3. Then we 
use a theorem of Hamilton \cite{Ha-et} to conclude that the 
final rescale limit is a steady soliton. 
Our arguments are similar to Perelman's  
compactness/convergence methods that were used extensively 
in his papers \cite{P1}, \cite{P2}; for reader's convenience we will 
give a detailed account. 

Recently it comes to our attention that Gu and Zhu \cite{GZ} 
proved the existence of type II singularity; they used Perelman's 
$l$ functional argument to detect type II singularity in the 
radial symmetric case. See also 
a very recent paper \cite{GI}.  

{\bf Acknowledgment}. 
We thank Prof. B. Chow, Prof. D. Knopf and Prof. G. Wei for 
helpful conversations. 
In June 2007 we gave a seminar talk at UCSB on this result, we are 
especially grateful to Prof. R. Ye for his questions and comments. 


\section{Notations and Definition}

All manifolds we consider in this paper are of dimension $3$. 
We use $R$, often with two variables $x$ and $t$, 
to denote the scalar curvature. $\Rm$ denotes the full curvature tensor. 
The Hamilton-Ivey  
pinching inequality (see \cite{Ha99}) says, 
if in the beginning the curvature 
is bounded from below by $-1$, then we have 
\[
\Rm \geq -\phi(|R|), 
\]
where $\phi$ is a nonnegative function so that 
$\lim_{r\rightarrow}\phi(r)/r=0$. In particular, when 
$R$ is large, the full curvature is dominated 
by the scalar curvature $R$.

We follow some notations of Perelman \cite{P1}, \cite{P2}. 
Assume $g(t)$ is a family of metrics on a manifold $M$  
that evolve under the Ricci flow.  
$B(x, t, r)$ denotes the metric ball centered at $x$, 
of radius $r$ with respect to the metric $g(t)$.  
Then one defines the {\em parabolic neighborhood} 
\begin{equation}
P(x, t, r, -\Delta t)= B(x, t, r)\times [t-\Delta t, t].
\end{equation} 
When we say two sets $(U_1, p_1, t_1)$ and $(U_2, p_2, t_2)$ 
are $\epsilon$-close, 
typically we first do rescales on $U_1$ and on $U_2$ so that 
$R(p_1, t_1)=R(p_2, t_2)=1$, then $\epsilon$-close means 
these rescaled sets are 
$\epsilon$-close under $C^{5}$ topology. There is a similar notion 
of $\epsilon$-close between two parabolic neighborhoods. 

Since 
we study only finite time singularities, by 7.3 of \cite{P1},  
all solutions are {\em noncollapsing}. In particular, when a sequence 
of parabolic neighborhoods admit a uniform curvature bound, we can 
take pointed limit over a subsequence. 

An {\em ancient solution} is a solution that 
exists on the time interval $(-\infty, T)$ for some $T\in {\mathbb R}$. 
The {\em cylinder} $S^2\times {\mathbb R}$, with the 
$S^2$ direction evolving under the Ricci flow,  is an important 
example. A set $(Z, z)$ evolving under the Ricci 
flow over time $[-t, 0]$ is called a {\em strong $\epsilon$-neck}, 
if after rescale the metrics 
by $R(z, 0)$, $Z$ is $\epsilon$-close to an evolving cylinder 
of length $\epsilon^{-1}$ from time $-\epsilon^{-1}$ to $0$.   

\begin{definition}
The {\bf caliber} of a cylinder $Y=S^2\times {\mathbb R}$ is  
$R^{-1}$, where $R$ is the scalar curvature of $Y$.
\label{defini-DS-caliber-cylinder-df93-gj912-sda8g}
\end{definition}
Clearly the caliber of a cylinder $Y$ is just 
half of the square of its radius. The following is obvious: 
\begin{lemma}
Start the Ricci flow at time $0$ 
on $Y$, then the caliber of $Y$ equals to the time it takes 
for the cylinder to go singular (i.e. 
shrink into the real line ${\mathbb R}$). 
\end{lemma} 
\begin{prop}[Perelman]
Given $\epsilon>0$. 
Assume the initial metric $g(0)$ satisfies the curvature 
bound $|\Rm (g(0))|\leq 1$, and for all $p\in M$, 
$\Vol_{g(0)}(B_1(p))\geq 10^{-1}$. 

Then there exists $r_0>0$, so that whenever $R(x_1, t_1)>r_0^{-2}$, 
the neighborhood $P(x_1, t_1, R(x, t)^{1/2}\epsilon^{-1}, 
-\epsilon^{-1}R(x, t))$, 
under the rescaled metric $R(x_1, t_1)g(t)$, 
is $\epsilon$-close to a parabolic neighborhood in one of the following: 

i). A space form with positive curvature evolving under the Ricci flow,  

ii). The cylinder $S^2\times {\mathbb R}$ (or 
$S^2\times {\mathbb R}/{\mathbb Z}_2$) 
evolving under the Ricci flow,

iii). A compact ancient solution 
with strictly positive, nonconstant (at each time slice) 
curvature that is diffeomorphic to $S^3$ or ${\mathbb {RP}}^3$,

iv). A noncompact ancient solution to the Ricci flow with 
strictly positive curvature. 
\label{prop-DS-canonical-ig93-9hg0-jka9-9j90-cjk3}
\end{prop}

For a proof, see theorem 12.1 of \cite{P1}, 
together with section 1 of \cite{P2}. 
The possibilities i), ii), iii), iv) above give 
a rough classification of noncollapsing 
ancient solutions of nonnegative curvature in dimension $3$. 
The parabolic neighborhood $P$ above is 
called a {\bf canonical neighborhood}. 

Both of the ancient solutions iii (in sufficiently ancient time), 
and iv (in all time), contain a piece of evolving cylinders; see 
\cite{P1} sections 11, 12 
and especially 1.4 of \cite{P2}. 
To emphasis the difference between 
cases ii), iii), iv), remember 
\begin{prop}[Perelman]
Assume $X$ is an 
ancient, noncollapsing $3$-dimensional ancient 
solution with nonnegative, 
nonconstant (at each time slice) curvature. Then at each time (e.g. $t=0$), 
$X$ can be decomposed into two parts, $X_C$ 
and $X_T$; either of them can be empty. 
$X_T$ is connected; $X_C$ has at most two connect 
components and is connected when $M$ is noncompact. The boundary 
components 
(if any) of $X_C$ and $X_T$ are all diffeomorphic to $S^2$. 

Each connected component of $X_C$ is compact, 
for all $p, q$ in the same connect component of $X_C$, we have 
\begin{equation}
R(p,0)\leq A(\epsilon)R(q, 0), \ \ \ \ 
\Diam X_C\leq D(\epsilon)R(p,0)^{-1/2}. 
\end{equation}
Every point in the second part, 
$X_T$, is the center of a strong $\epsilon$-neck. 

Moreover, when $X_C$ is connected, then it is diffeomorphic to 
either $S^3$, or $RP^3$, or the $3$-dimension ball $B^3$. When $X_C$ 
is not connected and therefore has two components, then one component 
of $X_C$ is diffeomorphic to $B^3$, the other is diffeomorphic to 
either $B^3$ or $RP^3-B^3$. 
\label{theo-DS-compact-tube-decomposition-gf9u39-fd89hf}
\end{prop}

Most of the above is proved in section 11 (especially 11.8) of \cite{P1}, 
the last part concerning the topology of $X_C$ can be found 
in section 1 of \cite{P2}; the proof uses a compactness 
argument (go to ancient time) involving the soul theorem, 
see for example chapter 9 of \cite{MT}, \cite{KL} and \cite{CaoZhu}. 

Roughly speaking, the above Proposition says the following.  
Assume $X$ is not a space form, then $X$ can be decomposed into a 
``tube part'' $X_T$ and a ``cap part'' $X_C$. 
More precisely, when $X$ is noncompact, 
then either $X$ is just the tube, or it is an approximate 
tube being connected, through one $S^2$ boundary, to a cap 
that is diffeomorphic to $B^3$ or $RP^3-B^3$. The later case is isometric 
to $S^2\times {\mathbb R}/{\mathbb Z}_2$. 

\begin{picture}(400,100)(0,0)
\put(80, 80){\line(1,0){270}}
\put(80, 50){\line(1,0){270}}
\qbezier(80,80)(30,80)(30,65)
\qbezier(80,50)(30,50)(30,65)
\qbezier(130,80)(120,65)(130,50)
\qbezier(130,80)(140,65)(130,50)
\put(220,60){$X_T$}
\put(80,60){$X_C$}
\put(180,10){\text{Figure I}.}
\end{picture}

\noindent 
Figure I shows the cap-tube decomposition
of the ancient solution case iv) in Proposition 
\ref{prop-DS-canonical-ig93-9hg0-jka9-9j90-cjk3}. 
We can take the convention that, 
if the canonical 
neighborhood in Proposition 
\ref{prop-DS-canonical-ig93-9hg0-jka9-9j90-cjk3} fall into case iv) 
{\em and} the canonical neighborhood lies in the $X_T$ part, then we shall 
classify this neighborhood into case ii). In another word, when the 
canonical neighborhood is of case iv), then the neighborhood contains 
$X_C$.

When $X$ is compact, 
we have two cases. If $X$ is not ``long enough'', then $X$ itself is 
$X_C$, which is diffeomorphic to $S^3$ or $RP^3$; if $X$ is ``long enough'', 
then $X$ comes from a bounded tube $X_T$ being capped at its two ends by 
two caps $X_C$: at one of these cap 
must be $B^3$, the other may be either $B^3$ or $RP^3-B^3$.  

\begin{picture}(400,100)(0,0)
\put(80, 80){\line(1,0){230}}
\put(80, 50){\line(1,0){230}}
\qbezier(80,80)(30,80)(30,65)
\qbezier(80,50)(30,50)(30,65)
\qbezier(310,80)(360,80)(360,65)
\qbezier(310,50)(360,50)(360,65)
\qbezier(280,80)(270,65)(280,50)
\qbezier(280,80)(290,65)(280,50)
\qbezier(110,80)(100,65)(110,50)
\qbezier(110,80)(120,65)(110,50)
\put(310,60){$X_C$}
\put(190,60){$X_T$}
\put(60,60){$X_C$}
\put(180,10){\text{Figure II}.}
\end{picture}

\noindent 
Figure II is a picture of the ancient solution iii in Proposition 
\ref{prop-DS-canonical-ig93-9hg0-jka9-9j90-cjk3}; this is 
the ``long'' case so that the tube part $X_T$ 
is nonempty.\footnote{
As we will see in the proof of Theorem 
\ref{th-DS-main-theorem-steady-soliton-faih-g9uj3-gve89}, 
case iii) can be ignored in our study.} 
In this case, if the canonical neighborhood falls on $X_T$, 
or, on a component of $X_C$ that is diffeomorphic to $RP^3-B^3$, 
then we classify it into case ii of Proposition 
\ref{prop-DS-canonical-ig93-9hg0-jka9-9j90-cjk3}; if it falls on a 
component of $X_C$ that is diffeomorphic to $B^3$, 
by a compactness argument we classify it into iv) of Proposition 
\ref{prop-DS-canonical-ig93-9hg0-jka9-9j90-cjk3}.

In particular, 
we can choose $\epsilon$ so that the options i), ii), iii), iv) 
of Proposition 
\ref{prop-DS-canonical-ig93-9hg0-jka9-9j90-cjk3} 
are {\bf mutually exclusive}. 
With these information, we give a working definition of 
{\em degenerate singularity} of Ricci flow in dimension $3$:

\begin{definition}
Assume $M$ is compact and 
$g(t)$ is a solution to the Ricci flow as in Proposition 
\ref{prop-DS-canonical-ig93-9hg0-jka9-9j90-cjk3} that exists on the 
time interval $[0, T)$ with $T<\infty$. 
Assume there is a sequence of points $(x_i, t_i)$ so that 
$\lim_{i\rightarrow\infty}t_i=T$ 
and $\lim_{i\rightarrow\infty}R(x_i, t_i)=\infty$; 
moreover for sufficiently large 
$i$ the canonical neighborhood of $(x_i, t_i)$ is of 
case iv) in Proposition \ref{prop-DS-canonical-ig93-9hg0-jka9-9j90-cjk3}.  
Then we say a degenerate singularity 
happens at time $T$. 
\label{Define-DS-degenerate-singularity-7834h-8as1}
\end{definition} 

A glance at Hamilton's picture in \cite{HaSing} suggests that one might 
also include case ii) 
in Proposition \ref{prop-DS-canonical-ig93-9hg0-jka9-9j90-cjk3}, when 
the canonical neighborhood is $S^2\times{\mathbb R}/{\mathbb Z}_2$; 
we will discuss this possibility in the end of this paper. 

It is possible that the solution goes 
singular {\em everywhere}, i.e. scalar curvature goes to 
infinity everywhere as $t\rightarrow T^-$. At this moment we don't 
know a compact example that extincts everywhere while developing 
a degenerate singularity. On the other hand, 
Perelman's {\em standard solution} in Section 2 of \cite{P2} is 
a noncompact example. 
If the solution extincts everywhere at $T$, 
then the topology of $M$ is quite simple, see sections 
3 and 4 of \cite{P2}.

There is an important gradient estimate for scalar curvature, see 
(1.3) of \cite{P2}: 
\begin{prop}[Perelman]
Let $g(t)$ be a solution to the Ricci flow as in 
Proposition \ref{prop-DS-canonical-ig93-9hg0-jka9-9j90-cjk3} and $r_0$ 
be the canonical neighborhood parameter. Then 
whenever $R(x, t)>r_0^{-2}$, we have 
\begin{equation}
|\nabla R(x, t)|\leq C_1R^{3/2}, \ \ \ \ 
|\partial_t R|(x, t)\leq C_2R^2.
\label{eq-DS-gradient-estimate-8bgg2-z19v-0lc0}
\end{equation}
\end{prop}

The next is a very useful locally splitting theorem:
\begin{prop}[Perelman] 
Given any $\epsilon$, 
there exists $\eta>0$ so that the following is true: 

Assume $X$ is a noncollapsing, noncompact, 
ancient solution of nonnegative, 
bounded curvature defined for time $t\in (-\infty, 0]$. 
Assume $\gamma$ is a minimal geodesic segment at time $0$, with 
two end points $p_1$ and $p_2$, and $p$ is a point on $\gamma$. 
Assume 
\begin{equation}
R(p_1, 0)>\eta^{-1}, \ \ \ \ \ \ \ \ 
R(p, 0)=1, \ \ \ \ \ \ \ \ 
R(p_2, 0)<\eta. 
\end{equation}

Then $d(p_1, p)>\epsilon^{-1}$, $d(p, p_2)>\epsilon^{-1}$, and 
$B(p, 0, \epsilon^{-1})$ is $\epsilon$-close to a subset in a 
cylinder of caliber $1$. 
\label{th-DS-splittin-theo-j83-f8jf8-27fjh34}
\end{prop}
This follows from Perelman's compactness theorem, see 11.7 of \cite{P1}: 
roughly, if this is not true, take a limit of counterexamples. By 
Perelman's curvature bound (see the last 
three lines in page 30 of \cite{P1}), we see the distance from $p$ to 
$p_1, p_2$ goes to $\infty$. Therefore the limit contains a line and  
splits; it must be the cylinder and that is a contradiction. 
See also \cite{KL}, \cite{MT}, \cite{CaoZhu} for more details.  

Finally, with nonnegative curvature, a singularity must happen everywhere: 
\begin{prop}[Perelman]
Assume $X$ is a noncollapsing 
ancient solution of strictly positive curvature. If $X$ develops 
a singularity at time $T$, then for all $x\in X$, 
\begin{equation}
\lim_{t\rightarrow T^-} R(x, t)=\infty.
\end{equation}
\end{prop}
See the proof of 12.1,  Claim 2, in \cite{P1}. 


\section{A theorem on ancient solutions} 

Very similar to Definition 
\ref{defini-DS-caliber-cylinder-df93-gj912-sda8g}, we define 
\begin{definition}
Assume $X$ is a noncollapsing, noncompact ancient 
solution of strict positive, bounded 
curvature, defined for time $(-\infty, 0]$. 
The {\bf caliber} of $X$ at time $0$ 
\begin{equation}
\lim_{x\rightarrow\infty}\frac{1}{R(x, 0)}.
\end{equation}
\end{definition}
In view of Proposition 
\ref{theo-DS-compact-tube-decomposition-gf9u39-fd89hf}, the above 
limit exists; in fact, the solution is close to a cylinder with 
the above caliber when $x\rightarrow\infty$. 

Since $(X, g(0))$ has bounded curvature, by Shi's theorem \cite{Shi1}  
we can extend the solution for a short time beyond $t=0$. 
Moreover, on each time slice, the curvature is bounded, positive;  
see \cite{Shi2} Theorem 4.14. See also 12.1 of 
\cite{MT} for an alternative argument. 

\begin{lemma}
Assume $g(t)$ is a maximal solution with initial data $(X, g(0))$ 
that has bounded curvature in each time 
slice. Then $g(t)$ exists on 
$[0, C)$, where $C$ is the caliber of $(X, g(0))$. Moreover, if 
$C<\infty$, then for all $x\in X$, 
$\lim_{t\rightarrow C^-}R(x, t)=\infty$. 
\label{lemma-DS-caliber-decides-life-ah92-sd8h-2k9f}
\end{lemma}
Here ``maximal'' means the following. We can give a partial order to the 
set of all Ricci flow solutions with bounded curvature at each time slice 
and with initial data $(X, g(0))$: we say $g_1\leq g_2$ if 
$g_2$ is an extension (in time) of $g_1$. Maximal means maximal according 
to this partial order. \footnote{
Implicitly in this statement, extension(s) may or may not 
be unique. Recent there have been results on 
the uniqueness of Ricci flow, these could make some of our arguments 
more simple.} 
\proof
The proof is exactly the same as the argument in section 2 of \cite{P2}, 
where Perelman proved the the life of the {\em standard solution} is its 
caliber at initial time. See \cite{KL}, \cite{MT}, \cite{CaoZhu} 
for more details. 
\qed

Recall that a solution develops a type I singularity, 
if the solution goes 
singular at time $T$ and there exists $C\geq 0$ so that 
\begin{equation}
\limsup_{t\rightarrow T^-} (T-t)\cdot \sup_{x\in X} R(x, t)=C.
\label{eq-DS-type-I-definition-neat-T-9fuqe-af893-gj8932}
\end{equation}

\begin{theorem}
There is no noncollapsing, noncompact ancient solution with strictly 
positive curvature that has a type I singularity. 
\label{th-DS-rule-out-type-II-ancient-sol-f983-fg943-9132wu}
\end{theorem}
We argue by contradiction. 
Assume there is one such solution $X$. Notice we assume positive 
curvature, so the solution is diffeomorphic to ${\mathbb R}^3$, its 
structure is described in Proposition 
\ref{theo-DS-compact-tube-decomposition-gf9u39-fd89hf}, see Figure I. 
In the next two lemmas we 
show that one can get a global curvature bound similar to 
(\ref{eq-DS-type-I-definition-neat-T-9fuqe-af893-gj8932}). 
\begin{lemma}
We have 
\begin{equation}
\liminf_{t\rightarrow T^-} (T-t)\cdot \sup_x R(x, t)>0.
\end{equation}
In particular, $C>0$ and there exists $c>0$ so that for 
all $t$ close to $T$ and all $x$,  
\begin{equation}
c\leq (T-t)R(x, t)\leq  C.
\end{equation}
\end{lemma}
\proof
By lemma \ref{lemma-DS-caliber-decides-life-ah92-sd8h-2k9f}, 
at time $t$ the caliber of $(X, g(t))$ is $T-t$. Therefore 
\begin{equation}
\lim_{x\rightarrow\infty} R(x, t)=\frac{1}{T-t}.
\end{equation}
We claim there is a universal constant $c$ so that for all $y\in X$, 
\begin{equation}
R(y, t) \geq c \lim_{x\rightarrow\infty} R(x, t)=\frac{c}{T-t}. 
\end{equation}
As in many occasions in Perelman's papers, e.g. the last paragraph 
in 11.7 of \cite{P1}, 
this follows from Yau's volume comparison argument with base point 
at infinity. For reader's 
convenience, we give a sketch:

By doing a rescale we can assume $T-t=1$. 
Take a ray 
$\gamma$ with $\gamma(0)=y$. If the distance $s_1$ is sufficiently 
large, a neighborhood of $\gamma(s_1)$ is close to a (piece of) 
cylinder with caliber $T-t$. Take $s_2\gg s_1$, let $p=\gamma(s_2)$, 
define the one-direction annulus 
\begin{equation}
A(p, \delta_1, \delta_2)
=\{x\in X, | \delta_1\leq d(x, p)\leq \delta_2, \,\, 
\liminf_{s\rightarrow \infty}d(\gamma(s), x)-(s-s_2)\geq 0 
\}. 
\end{equation}
The volume comparison based at $p$ implies 
\begin{equation}
\frac{\Vol(A(p, s_2-D, s_2+D))}{\Vol(A(p, s_2-s_1-D, s_2-s_1+D))}
\leq\frac{\int_{s_2-D}^{s_2+D} r^2 dr}{\int_{s_2-s_1-D}^{s_2-s_1+D} r^2 dr}
\approx 1, 
\label{eq-DS-vol-comp-9j0re3-g9hcf-gv0q9jt-g9}
\end{equation}
here $D\gg \sqrt{T-t}$ while $D\ll s_1\ll s_2$. Since 
$A(p, s_2-s_1-D, s_2-s_1+D)$, which contains $\gamma(s_1)$, is 
approximately a cylinder of length $2D$ and caliber $T-t$, we see 
\begin{equation}
\Vol(A(p, s_2-s_1-D, s_2-s_1+D))\sim 2D(T-t)\ll D^3. 
\label{eq-DS-vol-comp-asd90j4t-g90384h8hg84yr-2q09}
\end{equation}
If $R(y, t)$ is too small, say $R(y, t)\leq c$, 
then by Perelman's compactness theorem for 
ancient solution (see the proof of 11.7 of \cite{P1}), 
we have $R\leq 1$ on 
$B(y, t, D)$, so by the noncollapsing assumption we have 
\begin{equation}
\Vol(A(p, s_2-D, s_2+D))\geq 
\Vol(B(y, t, D))\geq \tau D^3.
\label{eq-DS-vol-comp-9ji4-g8hf-38t8u-fj930}
\end{equation}
Now (\ref{eq-DS-vol-comp-9ji4-g8hf-38t8u-fj930}) together with 
(\ref{eq-DS-vol-comp-asd90j4t-g90384h8hg84yr-2q09}) contradicts 
(\ref{eq-DS-vol-comp-9j0re3-g9hcf-gv0q9jt-g9}).
\qed

\begin{lemma}
If $X$ is of type I, then there exists a noncompact, 
noncollapsing ancient solution with 
positive curvature that goes singular at time $0$, 
and for all $x\in X_1$ and $t\leq 0$ we have 
\begin{equation}
c\leq |t|\cdot R(x, t)\leq C. 
\end{equation} 
\end{lemma}
\proof
Let $t_i\rightarrow T^-$. 
Take pointed rescale limit at points $(x_i, t_i)$ in the 
cap region, i.e. at the $X_C$ region 
(see Proposition \ref{theo-DS-compact-tube-decomposition-gf9u39-fd89hf}), 
this will guarantee the limit does not split. 
Then translate the singular time to $0$. 
The conclusion follows from the previous lemma. 
\qed

However, such ancient solution does not exist: 

\begin{theorem}
There is no noncollapsing, noncompact, 
ancient solution, with bounded positive curvature at 
each time slices, satisfying the following: 

1). The solution goes singular at time $T=0$. 

2). For all $(x, t)$ we have $c\leq |t|\cdot R(x, t)\leq C$. 
\label{th-DS-rule-out-type-one-ancient-noncpt-sol-afj4-h90j-sadf}
\end{theorem}

We break the proof into a sequence of lemmas. In view of  
Proposition \ref{theo-DS-compact-tube-decomposition-gf9u39-fd89hf}, 
at time $t$  when $x\in X_T$, we say $x$ {\bf is on the cap} at time $t$; 
when $x\in X_T$, we say $x$ {\bf is the center of a tube} at time $t$; 
see Figure I. When curvature is strictly 
positive, $X_C$ here is diffeomorphic to the solid ball $B^3$. 
A point $x$ may be on the cap at one time and be the center 
of a tube at another time. 
 
\begin{lemma}
Let $(p_i, t_i)$ be a sequence of points with $t_i\rightarrow -\infty$ 
that are on the cap. 
Then fix a time, say $-1$, the points $(p_i, -1)$ 
goes to space infinity when $i\rightarrow\infty$. 
\end{lemma}
\proof 
In fact, if this is not true, then 
assume $(p_i, -1)\rightarrow (p, -1)$.
Take $(p, -1)$ as the base point and take 
the constant curve $p$ from inverse time $0$ to $\tau$, 
and compute Perelman's $l$ functional, we get 
\begin{equation}
l(p, \tau)=\frac{1}{2\sqrt\tau}\int_0^\tau \sqrt s  R ds
\leq \frac{1}{2\sqrt\tau}\int_0^\tau \sqrt s  \frac{C}{s+c} ds
\leq C, \ \ \ \ \text{ for all } \ \ \tau>0.
\end{equation}
Thus by Perelman's asymptotic soliton theorem (\cite{P1} Proposition 11.2, 
see also \cite{KL}, \cite{MT}), 
if we take a rescaled limit at 
$(p, \tau)$, we shall get a nonflat shrinking soliton. The limit soliton  
cannot be compact, and since $p$ stay on the cap all the time,  
the limit does not split, i.e. the limit soliton 
is of positive curvature. By section 1 of \cite{P2}, 
shrinking solitons are either space forms or 
$S^2\times {\mathbb R}$, or $S^2\times {\mathbb R}/ {\mathbb Z}_2$,  
none of them contains a cap-like neighborhood. That is a contradiction. 
\qed

\begin{lemma}
Given any $C>0$, there exists $\epsilon$ so that for any $\tau>2$ and 
for any ancient solutions $(M, p)$ with $\sup R\leq C$ at time $-1$,  
if 
\begin{equation}
\tilde V(\tau)\leq \tilde V(2\tau)+\epsilon,
\end{equation}
and $l(q, \tau)\leq C$, with the base point of $l$ and $\tilde V$ 
taken at $(p, -1)$, 
then after rescale by $R(q, \tau)$, a neighborhood of $(q, \tau)$ 
is close to a subset in 
a shrinking soliton. Here $\tau=-1-t$ is the inverse time. 
\label{lemma-DS-almost-max-reduced-volume-9ufv-fv8j-bl23}
\end{lemma}
\proof
This follows exactly the proof for 
Proposition 11.2 in \cite{P1}. We just need 
to replace the bound $n/2$ for $l$ there by $C$. See \cite{Ye} 
for more details. 
\qed

\begin{lemma}
Assume there is an ancient solution $X$ as in 
Theorem \ref{th-DS-rule-out-type-one-ancient-noncpt-sol-afj4-h90j-sadf}. 
Then for any $\epsilon>0$, there exists 
$N>0$, so that  
for any point $(p, -1)$ that is sufficiently far into space infinity, 
we have 
\begin{equation}
\tilde V(\tau)\leq \tilde V (2\tau)+\epsilon, 
\end{equation}
for all $\tau>N$. 
Here we take $(p, -1)$ as the base point for computing 
$l$ and $\tilde V$.
\label{eq-SD-almost-max-V-90kh98-g84hd-238fhhe8j3954}
\end{lemma}
\proof 
Observe, at time $t$, at space infinity the solution 
looks like a tube of caliber $|t|$.   
So no matter where the base point $(p, t)\in X$ for 
$l$ locates, we have 
\begin{equation}
\tilde V_X(\infty)=\lim_{\tau\rightarrow\infty}\tilde V_X(\tau)
=\frac 2e(4\pi)^{3/2}. 
\label{eq-DS-ancient-tilde-V-value-cylinder-asdf-349ujgtf-gfw0}
\end{equation}
In fact, pick a base point $(p, t)$, 
translate the time $t$ to $0$; then we have the inverse time $\tau$. 
Apply Proposition 11.2 in \cite{P2}, because 
the asymptotic soliton is noncompact, it must be $S^2\times {\mathbb R}$ 
(it is easy to check $S^2\times {\mathbb R}/{\mathbb Z}_2$ 
does not happen). Therefore the ancient soliton is a evolving cylinder 
that goes singular at time $0$. On the cylinder we see 
$\tilde V_{\text{cylinder}}(\tau)=\tilde V_{\text{cylinder}}(\infty)$ 
does not depend on $\tau$. Then 
because $\Ric+\Hess_l={g}/{(2\tau)}$ on the cylinder,  
it is easy to compute that when $\tau=1$, 
\begin{equation}
l_{\text{cylinder}}=1+x^2/4, 
\end{equation} 
where $x$ is the coordinate in the ${\mathbb R}$ direction. 
Therefore integrate on the cylinder, 
\begin{equation}
\tilde V_{\text{cylinder}}(\infty)=\tilde V_{\text{cylinder}}(1)
=\int_{\text{cylinder}} e^{-l}=\frac 2e(4\pi)^{3/2}.
\end{equation}
By Perelman's growth estimate for $l$ (see Lemma 3.2 in \cite{Ye}), 
one can bound the contribution to 
$\tilde V$ outside any give distance. Moreover the scale-down 
limit of $l$ is just $l_{\text{cylinder}}$. 
Then one gets
\begin{equation}
\tilde V_X(\infty)
=\tilde V_{\text{cylinder}}(\infty)=2e^{-1}(4\pi)^{3/2}.
\end{equation}

Especially, 
pick any $(z, -1)$ in an evolving cylinder $Z$ 
as the base point for computing $l$, we get the same limit value 
$\tilde V_Z(\infty)=2e^{-1}(4\pi)^{3/2}$. 
Although both $\tilde V_Z$ and 
$\tilde V_{\text{cylinder}}$ are computed on evolving cylinder, 
they are different in that $\tilde V_{\text{cylinder}}$ is defined as 
the limit of reduced volumes over rescaled manifolds, essentially it 
used the singularity as base point; 
$\tilde V_Z$ used an ordinary point $(z, -1)$ 
as base point. Therefore $\tilde V_{\text{cylinder}}(\tau)$ does not 
depend on $\tau$ while $\tilde V_Z(\tau)$ does. 
 
There exists $N$ so that 
\begin{equation}
\tilde V_Z(N)\leq 2e^{-1}(4\pi)^{3/2}+\epsilon/3. 
\label{eq-DS-cylinder-84q23-g8j9q-sa8934}
\end{equation}
Moreover,  there is a radius $D\gg \sqrt N$ so that 
\begin{equation}
\tilde V_Z(N)-\frac{\epsilon}{3}\leq \int_{B(z, -1-N, D)} 
N^{-n/2}e^{-l}\leq 
\tilde V_Z(N)\leq \frac 2e(4\pi)^{3/2}+\frac {\epsilon}{3}. 
\label{eq-DS-cylinder-9iu32-gh73aq2d-f89jsd}
\end{equation}

If we take pointed limit with base points $(p_i, -1)\in X$ 
going to space infinity, 
we get a cylinder of caliber $1$.  
Therefore if we go sufficiently far into the space infinity, 
we see the ball $B$ centered at $(p, -1)$ with radius $D$ satisfies that 
$B\times [-1-N, -1]$ is close to the corresponding part in space-time 
of the evolving cylinder $Z$. So 
by Perelman's growth estimate on $l$ (see Lemma 3.2 of 
\cite{Ye}), we have 
\begin{equation}
\tilde V_X(N)-\frac{\epsilon}{3}\leq 
\int_{B(p, -N-1, D)} N^{-n/2}e^{-l_X} < \tilde V_X(N), 
\end{equation}
\begin{equation}
 \Big|\int_{B(p, -N-1, D)}  N^{-n/2}e^{-l_X}
- \int_{B(z, -N-1, D)}  N^{-n/2}e^{-l_Z}\Big|\leq \frac{\epsilon}{3}. 
\end{equation}
Recall $\tilde V$ is a decreasing function in $\tau$; 
see (7.12)-(7.13) in \cite{P1}. So the above combined with 
(\ref{eq-DS-ancient-tilde-V-value-cylinder-asdf-349ujgtf-gfw0}), 
(\ref{eq-DS-cylinder-84q23-g8j9q-sa8934}) 
and 
(\ref{eq-DS-cylinder-9iu32-gh73aq2d-f89jsd}) 
proves the conclusion. 
\qed

Now we prove Theorem  
\ref{th-DS-rule-out-type-one-ancient-noncpt-sol-afj4-h90j-sadf}. 
Take a sequence $(p_j, t_j)$ with $t_j\rightarrow -\infty$ that 
are on the cap.  
Then for sufficiently large $j$ with $|t_j|>N$ (here $N$ is the number in 
Lemma \ref{eq-SD-almost-max-V-90kh98-g84hd-238fhhe8j3954}), 
we see $(p_j, -1)$ is so far 
into space infinity that if we use $(p_j, -1)$ as the bast point for $l$, 
we get 
\begin{equation}
\tilde V(\tau)\leq \tilde V (2\tau)+\epsilon;
\end{equation}
for all $\tau>N$. In particular we take $(p_j, -1)$ as the base point 
and use $\tau=-1-t_j$. The type I assumption implies a bound of $l$ 
at $(p_j, t_j)$, evaluated on the constant curve $p_j$.     
Therefore by Lemma 
\ref{lemma-DS-almost-max-reduced-volume-9ufv-fv8j-bl23}, after rescale 
by $R(p_i, t_i)$ we see a big neighborhood of $(p_i, t_i)$ is close to 
s set in a shrinking soliton. But this is impossible because 
$(x_i, t_i)\in X_C$.  
\qed

%
%
%
%
%
%


\section{Forward limit and eternal solutions}
To get more concrete 
information, especially for taking forward limit, 
we need the following lemma, which could be used as an 
alternative definition of degenerate singularity; this is very 
similar to Lemma 4.3 in Perelman's second paper \cite{P2}. 

\begin{lemma}
Assume a degenerate singularity happens at time $T$, 
and the solution does not extinct everywhere at $T$. 
Then there exists a positive real number $\epsilon>0$ (depending on 
the solution) and a compact set 
$B\subset M$ that is diffeomorphic to a solid $3$-ball, so that 
there exists 
a neighborhood of $\partial B$, for all $t\in [T-\epsilon, T)$, 
after rescale the metric by $R(p, t)$, is $10^{-2}$-close to a 
cylinder of scalar curvature $1$ and length $100$; and we have 
\begin{equation}
\limsup_{t\rightarrow T^-}\sup_{p\in \partial B} R(p, t)<\infty. 
\end{equation}
Moreover, let $p(t)\in B$ be a point so that 
\begin{equation}
d_{g(t)}(p(t), \partial B)=\sup_{p\in B} 
d_{g(t)}(p, \partial B),
\label{eq-DS-find-the-tip-of-cap-dsfh9-g899e4-fj9843}
\end{equation}
then 
\begin{equation}
\lim_{t\rightarrow T^-}R(p(t), t)=\infty.
\end{equation}
\label{lemma-DS-alt-define-degenerate-sing-h4g3-fg9j3-fh873ed3}
\end{lemma}
\proof
Let $(x_i, t_i)$ be a sequence of points as in Definition 
\ref{Define-DS-degenerate-singularity-7834h-8as1}. Since $M$ is compact, 
we can assume $x_i\rightarrow x_\infty$. Observe 
\begin{equation}
\lim_{t\rightarrow T^-} R(x_\infty, t)=\infty.
\end{equation}
In fact, if this is not true, then for some $C>0$, there exists $t$ 
so that $T-t$ is arbitrarily small while $R(x, t)<C$; we pick such a 
time $t^*$ that is sufficiently close to $T$. 
By the gradient estimate  
(\ref{eq-DS-gradient-estimate-8bgg2-z19v-0lc0}) in the space direction, 
we see a neighborhood $\mathcal N$ 
of $(x_\infty, t^*)$ has bounded curvature. Now because we have 
chosen $t^*$ sufficiently close to $T$, by the time-direction 
gradient estimate in 
(\ref{eq-DS-gradient-estimate-8bgg2-z19v-0lc0}), $\mathcal N$ has 
a uniform curvature bound on the time interval $(t^*, T)$. So 
on $\mathcal N$, the metrics from time $t^*$ to $T$ are uniformly 
equivalent. In particular, for sufficiently large $i$ we have 
$x_i\in \mathcal N$; this contradicts to the fact 
$R(x_i, t_i)\rightarrow\infty$. 

It is assumed that the solution does not extinct everywhere 
at time $T$. As we have seen above, by the gradient estimate 
(\ref{eq-DS-gradient-estimate-8bgg2-z19v-0lc0}), near time $T$,   
curvature cannot increase/decrease suddenly. So  
there is a point $y$ so that 
\begin{equation}
\limsup_{t\rightarrow T^-} R(y, t)=C<\infty.
\end{equation}
Pick a time $t^*$ that is sufficiently to $T$ so that 
\begin{equation}
R(x_\infty, t^*)>\eta^3\cdot \max\{ r_0^{-2}, C\},
\end{equation} 
here $r_0$ is the canonical neighborhood parameter in Proposition 
\ref{prop-DS-canonical-ig93-9hg0-jka9-9j90-cjk3}, and $\eta$ is the 
constant in Perelman's splitting argument, see Proposition 
\ref{th-DS-splittin-theo-j83-f8jf8-27fjh34}. 
Take a minimal geodesic $\gamma$ 
with respect to $g(t^*)$, from $x_\infty$ 
to $y$, parametrized by $s$. 
Let $r_1$ be so that $r_1^{-2}=\eta\cdot \max\{ r_0^{-2}, C\}$, and 
\begin{equation}
s_0=\inf \{s | R(\gamma(s), t^*)=r_1^{-2}\}, \ \ \ \ 
s_1=\inf \{s | R(\gamma(s), t^*)=4 r_1^{-2}\}.
\end{equation}
Now since $t^*$ 
is sufficiently close to $T$, by applying the gradient estimate 
(\ref{eq-DS-gradient-estimate-8bgg2-z19v-0lc0}), we can assume on the 
time interval $(t^*, T)$ we have 
\begin{equation}
R(\gamma(s_0), t)\leq 2r_1^{-2}.
\end{equation}

By Proposition \ref{th-DS-splittin-theo-j83-f8jf8-27fjh34}, 
the canonical neighborhood of 
$\gamma(s_0)$ at time $t^*$ is a tube, with radius $r_1\sqrt 2$. 
We now remove a center sphere $S^2$ in this tube at $\gamma(s_0)$ from 
the manifold $(M, g(t^*))$ and get an incomplete 
Riemannian manifold, possibly with $2$ components. 
In any case, pick the component that contains $x_\infty$ and take 
its closure $M^*$. 
$\partial M^*$ is either connected or has two components, all 
diffeomorphic to $S^2$. Write $S=\partial M^*$ 
when $\partial M^*$ is connected;  when $\partial M^*$ has 
two components, 
let $S$ the component of $\partial M^*$ 
that is closer to $\gamma(s_0-\delta)$, with 
$\delta= R(\gamma(s_0), t^*)^{-1/2}$. In another word, $S$ is 
directly hit by the part of $\gamma$ coming from $x_\infty$. 
 
The Ricci flow still runs on this manifold (it runs 
on the original manifold $M$, we just remove a sphere).  
For all $t\in [t^*, T)$, pick a point $q(t)$ that is furthest from 
$S$: 
\begin{equation}
d_t(q(t), S)=\max \{ d(x, S)\,|\, x\in M^*\}.
\end{equation} 
Let $\alpha$ be a minimum $g(t)$-geodesic from $S^2$ to $q(t)$. 
At time $t^*$, $S$ is the start of a piece of tube, 
in which $\gamma$ is passing through, with length $|s_1-s_0|$, 
along which the scalar curvature changes from 
$r_1^{-2}$ to $4r_1^{-2}$.  
There are two cases:

{\bf Case I}. 
For all $s$ with 
\[
s\geq |s_1-s_0|,
\]
we have 
\begin{equation}
R(\alpha (s), t^*)\geq 2 r_1^{-2}. 
\end{equation}
This means, if $\partial M^*$ has another component $S'$, then the 
canonical neighborhoods of $\alpha(s)$ does not intersect $S'$; in fact 
a neighborhood of $S'$ is very close to a long tube with scalar curvature  
$r_1^{-2}$. By the classification of canonical neighborhood and Proposition 
\ref{th-DS-splittin-theo-j83-f8jf8-27fjh34}, $M^*$ is being covered 
by canonical neighborhoods and therefore $\partial M^*=S$ is connected. 
Notice the neighborhood near $q(t^*)$ is necessarily a cap. 

In particular, 
$M^*$ is diffeomorphic to a solid ball or $RP^3-B$. 
Now $x_\infty\in M^*$, and the curvature at $\partial M^*=S^2$ remain 
bounded when $t\rightarrow T$, we conclude that for $i$ sufficiently 
large, $x_i\in M^*$. In particular, 
the possibility $M^*=RP^3-B$ is ruled out, because $M^*$ contains a 
cap diffeomorphic to $B^3$.  
Therefore, we can just make $p(t)=q(t)$ and $B=M^*$. 

{\bf Case II}. 

\begin{picture}(400,100)(0,0)

\qbezier(40,80)(120,75)(200,70)
\qbezier(40,40)(120,45)(200,50)
\qbezier(360,80)(280,75)(200,70)
\qbezier(360,40)(280,45)(200,50)
\qbezier(360,40)(370,60)(360,80)
\qbezier(360,40)(350,60)(360,80)
\qbezier(40,60)(220,50)(380,70)
\put(370,45){$S$}
\put(350,90){$\gamma(s_0)$}
\put(200,80){$\gamma(s_1)$}
\put(10,60){$\alpha(s_2)$}
\put(190,30){$R\approx 4r_1^{-2}$}
\put(350,25){$R\approx r_1^{-2}$}
\put(25,25){$R\approx 3r_1^{-2}$}
\put(200,80){$\gamma(s_1)$}
\put(180,10){\text{Figure III}.}
\end{picture}

There is a minimal $s'$ with 
$s>|s_0-s_1|$ so that 
\begin{equation}
R(\alpha(s'), t^*)=2r_1^{-2}.
\end{equation}
We will prove that this is impossible. Assume it does happen, 
we will find a minimum $s_2$ with $s_2>|s_0-s_1|$ so that 
\begin{equation}
R(\alpha(s_2), t^*)=3r_1^{-2}.
\end{equation}
Since we choose $t^*$ to be sufficiently close to $T$, 
by the gradient estimate 
(\ref{eq-DS-gradient-estimate-8bgg2-z19v-0lc0}), a neighborhood 
$\mathcal N_2$ of $\alpha(s_2)$ has a curvature bound 
\begin{equation}
R\leq 10r_1^{-2} 
\end{equation}
for all time $t\in [t^*, T)$.
Again by Proposition \ref{th-DS-splittin-theo-j83-f8jf8-27fjh34}, 
the canonical neighborhood near $\mathcal N_2$ is a cylinder. 
Now we cut $M^*$ at a center sphere at $\alpha(s_2)$ and take the component 
$\mathcal C$ that contains $\gamma(s_1)$.
This component $\mathcal C$ 
is just a piece of a cylinder, see also 
Proposition \ref{th-DS-splittin-theo-j83-f8jf8-27fjh34}. In 
particular, by the choice of $\gamma(s_1)$,  
we have $x_\infty\in \mathcal C$.  
Now we have uniform curvature bounds on the 
two ends of this cylinder on the time interval $[t^*, T)$. 
In particular, for sufficient large $i$, we have $x_i\in {\mathcal C}$. 
On the other hand, also by (\ref{eq-DS-gradient-estimate-8bgg2-z19v-0lc0}), 
there is actually a lower bound 
\begin{equation}
\frac {1}{10}r_1^{-2}\leq R(x, t) 
\end{equation}
for all $x\in{\mathcal C}$ and $t\in [t^*, T)$. This means for all time 
$t\in [t^*, T)$, all $x\in {\mathcal C}$ has a canonical neighborhood. 
We already know that near the two ends, the canonical neighborhoods 
are both cylinders, therefore we conclude all points in ${\mathcal C}$ 
has cylinder as a canonical neighborhood. This contradicts the fact that 
the canonical neighborhood of $(x_i, t_i)$ is a cap. 
\qed

Assume a degenerate singularity happens at time $T$. 
By taking a rescale limit at $(p(t), t)$ above, 
with $t\rightarrow T^-$,  
we get an ancient solution $X$. Since 
the canonical neighborhood of $p(t)$ are of case iv) in Proposition 
\ref{prop-DS-canonical-ig93-9hg0-jka9-9j90-cjk3}, 
the limit $X$  
is necessarily of strictly positive curvature, otherwise the solution is 
$S^2\times {\mathbb R}$ or its ${\mathbb Z}_2$ quotient, by the strong 
maximum principle. By a translation in time, $X$ 
exists on $(-\infty, 0]$. 

Let $T_1>0$ be the maximum time so that there is an 
extension of $X$ up 
to time $T_1$ and for all $t<T_1$, the curvature of $g(t)$ is positive 
and bounded. Therefore $T_1$ is the caliber of $(X, g(0))$, see 
Lemma \ref{lemma-DS-caliber-decides-life-ah92-sd8h-2k9f}. 
We are not saying this extension is unique. Moreover, 
{\em a priori}, when $0<T'<T_1$ we do not know if the solution 
$(X, g(t))$, with $-\infty< t\leq T'$, is a rescale limit of 
the original Ricci flow on $M$. 

\begin{lemma}
There exists an extension $(X, g(t))$ with $-\infty< t <T_1$, where $T_1$ 
is the caliber of $(X, g(0))$, so that 
for all $T'>0$ with $T'<T_1$, the ancient solution $X$ over time 
$(-\infty, T']$ is a rescale limit of the original solution on $M$. 
\label{lemma-SD-resolution-of-ancient-sol-gfr9e-32r57ed-gj923}
\end{lemma}
\proof 
Consider the set of time $T$, so that there is an extension 
of $X$ to time $T$, and the ancient solution $X$ over time 
$(-\infty, T)$ is a rescale limit of the original solution on $M$.
Let $T_0$ be the supremum of this set. 

First of all, we have $T_0>0$. In fact, a rescale of $M$ is close 
to $(X, g(0))$ under the pointed $C^3$ topology. 
We can assume the rescale is taken at $p(t)$ 
in Lemma 
\ref{lemma-DS-alt-define-degenerate-sing-h4g3-fg9j3-fh873ed3}, 
remember that the original solution does not go singular unless 
$R(p(t), t)$ goes to infinity. Now the corresponding neighborhood 
$\mathcal N(t)$ of $p(t)$ 
contains a cap part ${\mathcal N}_C$ and 
a tube part $\mathcal N_T$ that corresponds to the 
$X_C, X_T$ decomposition of $X$. 

From the definition of $p(t)$ in 
(\ref{eq-DS-find-the-tip-of-cap-dsfh9-g899e4-fj9843}), as long 
as the Ricci flow is running (at least for a short time when 
the metric geometry of ${\mathcal N}$ does not change too much), 
$p(t+\Delta t)$, 
the global maximum point of a distance function, must remain in 
$\mathcal N$. The 
gradient estimate (\ref{eq-DS-gradient-estimate-8bgg2-z19v-0lc0}) 
implies the curvature on $\mathcal N$ cannot go to infinity immediately; 
so $R(p(t+\Delta t), t+\Delta t)$ does 
not go to infinity immediately. During this time period 
the solution cannot go singular. Therefore we can 
take a {\em definite} 
short time forward limit at $p(t)$. This limit extends $X$, 
so $T_0>0$.

Now we can repeat this argument (i.e. view $T_0$ as time ``$0$'') as 
long as $X$ does not go singular, therefore we must have $T_0=T_1$, 
the caliber of $(X, g(0))$. 
\qed

We would like to have $T_1=\infty$, that is, $X$ is an eternal solution 
by Lemma \ref{lemma-DS-caliber-decides-life-ah92-sd8h-2k9f}. 
Assume this is not the case, then 
it is necessary that the solution develops a singularity at time $T_1$. 
By Theorem \ref{th-DS-rule-out-type-II-ancient-sol-f983-fg943-9132wu}, 
this  singularity is of type II as defined 
in \cite{HaSing}, that is, 
\begin{equation}
\limsup_{t\rightarrow T^-} (T_1-t)\cdot \sup_x R(x, t)
=\infty.
\end{equation}
In all cases, 
by taking a further limit we will get an eternal solution in which 
scalar curvature reaches maximum. This 
{\em point-picking} procedure is well known; \cite{P1}, \cite{P2}, 
\cite{KL}, \cite{CLN}. We 
give a quick account here.

We will find a sequence of points $(p_n, t_n)$, so that when 
we rescale the solution by $R(x_n, t_n)$ (i.e. first 
make $R(x_n, t_n)=1$), and then translate the time $t_n$ into $0$ 
(i.e. our base point is now $(p_n, 0)$ with $R(p_n, 0)=1$), 
then the solution exists on the interval $(-\infty, n]$, moreover, 
\begin{equation}
1=R(p_n, 0)\geq \Big(1-\frac{1}{n+1}\Big)
\sup\{R(x, t)\,|\, x\in X, t\leq n\}.
\end{equation}
In fact, 
start with $(p, 0)=(p^{(0)}, t^{(0)})$ 
so that $R(p, 0)$ is maximal at time $0$, and 
\begin{equation} 
(T_1- t^{(0)})R(p^{(0)}, t^{(0)})>(n+1)^2
\end{equation}
if $X$ is not eternal. Assume 
there is a point $p^{(1)}$ and time $t^{(1)}$ so that 
\begin{equation}
R(p^{(1)}, t^{(1)})>\Big(1+\frac 1n\Big)R(p^{(0)}, t^{(0)}), \ \ \ \ \ \ \ \ 
t^{(1)}-t^{(0)}\leq nR(p^{(0)}, t^{(0)})^{-1}.
\end{equation}
Without lose of generality we assume $R(p^{(1)}, t^{(1)})$ is 
almost maximal 
at the time slice $t=t^{(1)}$. If $(p^{(1)}, t^{(1)})$ does not satisfy 
our requirement, we continue to find $(p^{(2)}, t^{(2)})$...
In general, 
\begin{equation}
R(p^{(k)}, t^{(k)})>\Big(1+\frac 1n\Big)^kR(p^{(0)}, t^{(0)}), 
\ \ \ \ \,\,\,\,
t^{(k)}-t^{(0)}\leq n(n+1)R(p^{(0)}, t^{(0)})^{-1}.
\end{equation}
This is a contradiction, since the solution $X$ exists up to 
time $n(n+1)R(p^{(0)}, t^{(0)})^{-1}$, therefore there is a 
curvature bound at the time slice 
$t=n(n+1)R(p^{(0)}, t^{(0)})^{-1}$ 
(see also Theorem 11.4 in \cite{P1}). By Hamilton's 
Harnack inequality \cite{Ha-ha}, 
this scalar curvature bound actually holds for all 
$t\leq n(n-1)R(p^{(0)}, t^{(0)})^{-1}$.

Now take a pointed limit $Z$ of $(X, (p_n, t_n), R(p_n, t_n)g)$, we get 
an eternal solution on which the maximal of scalar curvature 
is reached. Therefore we can apply Hamilton's theorem 
on eternal solutions (see \cite{Ha-et}) 
to conclude that the limit $Z$ is a steady gradient soliton. Now by 
Lemma \ref{lemma-SD-resolution-of-ancient-sol-gfr9e-32r57ed-gj923} 
we finally get 

\begin{theorem}
Assume $g(t)$ develops a degenerate singularity at time 
$T$. Then 
then a rescale limit of $g(t)$ towards time $T$ 
is a steady gradient soliton. 
\label{th-DS-main-theorem-steady-soliton-faih-g9uj3-gve89}
\end{theorem} 
\proof
We already proved the first case, i.e. the solution 
does not extinct everywhere at time $T$. 

The second case, that $M$ extincts 
everywhere at time $T$, is actually easier. In fact the solution 
does not extinct 
unless the everywhere the curvature goes to infinity; therefore 
there is no difficulty in obtaining forward time limit. On the other hand, 
when we take limit as before and get an ancient 
solution $X$ with nonnegative curvature, then $X$ is 
noncompact. In fact if this is not true, then $X$ is diffeomorphic to 
$S^3$ or $PR^3$; by taking a forward limit we see that  
the {\em original} solution must be of strictly positive curvature near 
time $T$, in particular it will get rounder and rounder and cannot be 
a degenerated singularity. So the proof goes like the first case. 
\qed

\begin{corollary}
A rescale limit of Perelman's standard solution 
is the Bryant soliton. 
\end{corollary}
\proof
This follows from Theorem 
\ref{th-DS-rule-out-type-II-ancient-sol-f983-fg943-9132wu}. 
Notice that the soliton we get is also radial symmetric; see 
\cite{P2}, chapter 12 of \cite{MT}. 
Therefore by the discussion in \cite{CLN} 
the soliton must be the Bryant soliton.
\qed

Perelman conjectured that Bryant soliton is the only noncompact 
ancient solution with positive curvature; this will imply 
all the results in our paper. 
It should be true that every rescale limit taken along $(p(t), t)$ 
in Lemma 
\ref{lemma-DS-alt-define-degenerate-sing-h4g3-fg9j3-fh873ed3} 
must be the Bryant soliton. Unfortunately currently there are 
several difficulties in this direction. 

Finally we mention that analogue result of Lemma 
\ref{lemma-DS-alt-define-degenerate-sing-h4g3-fg9j3-fh873ed3} 
for $RP^3-B^3$ caps  also holds. In fact one can prove that locally 
the singularity is of type I and a rescale limit is 
$S^2\times {\mathbb R}/{\mathbb Z}_2$; one can get an example for 
such a singularity (e.g.  
taking ${\mathbb Z}_2$ quotient of the examples in chapter 2, \cite{CK}). 
Therefore we will not classify this case as degenerate.



\begin{thebibliography}{Al}

\bibitem{AK1} S. Angenent,D. Knopf, \textit{
An example of neckpinching for Ricci flow on $S\sp {n+1}$}.  
Math. Res. Lett.  11  (2004),  no. 4, 493--518.
\bibitem{AK2} S. Angenent,D. Knopf, \textit{
Precise asymptotics of the Ricci flow neckpinch}. 
arXiv:math/0511247v1.
\bibitem{CaoZhu} 
H. Cao, X. Zhu, 
\textit{Hamilton-Perelman's Proof 
of the Poincaré Conjecture and the Geometrization Conjecture}  
Asian J. Math., 10, (2) (2006), 165--492. math.DG/0612069.
\bibitem{CK} B. Chow, D. Knopf, \textit{Ricci flow, an introduction}
American Mathematical Society, 2004. 
\bibitem{CLN} B. Chow, P. Lu, L. Ni, \textit{Hamilton's Ricci flow}.
American Mathematical Society and Science Press, 2006.  
\bibitem{GI} D. Garfinkle, J. Isenberg, \textit{
The Modelling of Degenerate Neck Pinch Singularities 
in Ricci Flow by Bryant Solitons}. 
arXiv:0709.0514.
\bibitem{GZ} H. Gu, X. Zhu, \textit{
The Existence of Type II Singularities for the Ricci Flow on $S^{n+1}$}.
arxiv.0707.0033v1.
\bibitem{Ha-ha} R. Hamilton, \textit{The Harnack estimate for the 
Ricci flow}. J. Diff. Geom. 37 (1993) 225-243.
\bibitem{Ha-et} R. Hamilton, \textit{Eternal 
solutions of the Ricci flow}. J. Diff. Geom. 38 (1993) 1-11.
\bibitem{HaSing} R. Hamilton, \textit{The formation of singularities
in the Ricci flow}. Surveys in Differential Geometry II.
\bibitem{Ha99} R. Hamilton, \textit{Nonsingular solutions of
the Ricci flow on three-manifolds}. Comm. Anal. Geom. 1999. 
\bibitem{KL} B. Kleiner, J. Lott,
\textit{Notes on Perelman's papers}. math.DG/0605667. 
\bibitem{MT} J. Morgan, G. Tian,  
\textit{Ricci flow and the Poincar\'e conjecture}. 
Clay Mathematics Monographs, 2007. 
\bibitem{P1} G. Perelman,
\textit{The entropy formula for the Ricci flow and its geometric
applications}. math.DG/0211159. 
\bibitem{P2} G. Perelman, 
\textit{Ricci flow with surgery on three-manifolds}. math.DG/0303109.
\bibitem{Shi1} W. X. Shi, \textit{Deforming the metric on
complete Riemannian manifolds}. J. Diff. Geom. 30 (1989) 223-301.
\bibitem{Shi2} W. X. Shi, \textit{Ricci Deformation of the
metric on complete noncompact Riemannian manifolds}.
J. Diff. Geom. 30 (1989) 303-394.
\bibitem{Ye} R. Ye, \textit{On the $l$ function and the 
reduced volume of Perelman, I}. Download at 
http://www.math.ucsb.edu/$\sim$yer
\end{thebibliography}
\end{document}